\documentclass[12pt]{article}%
\usepackage{amsmath}
\usepackage{amsfonts}
\usepackage{amssymb}
\usepackage{graphicx}%
\providecommand{\U}[1]{\protect\rule{.1in}{.1in}}

\begin{document}
\title{\bf Analysis on a fractal set }

\author{
Santanu Raut and Dhurjati Prasad Datta\thanks{Corresponding author, email:dp${_-}$datta@yahoo.com} \\
Department of Mathematics, 
University of North Bengal \\
Siliguri- 734013, India } 
\date{}
\maketitle

\baselineskip =19pt

\begin{abstract}
The formulation of a new analysis on a zero measure Cantor set $C ( \subset I=[0,1])$  is presented. A non-archimedean absolute value is introduced in $C$ exploiting the concept of {\em relative} infinitesimals and a scale invariant ultrametric valuation of the form $\log_{\varepsilon^{-1}} (\varepsilon/x) $ for a given scale $\varepsilon>0$ and infinitesimals $0<x<\varepsilon, \ x\in I\backslash C$. Using this new absolute value, a valued (metric) measure is defined on $C $ and is shown to be equal to the finite Hausdorff measure of the set, if it exists. 
The formulation of a scale invariant real analysis is also outlined, when the singleton $\{0\}$ of the real line $R$ is replaced by a zero measure Cantor set. The Cantor function is realised as a locally constant function in this setting. The ordinary derivative $dx/dt$ in $R$ is replaced by the scale invariant logarithmic derivative $d\log x/d\log t$ on the set of valued infinitesimals. As a result, the ordinary real valued functions are expected to enjoy some novel asymptotic properties, which might have important applications in number theory and in other areas of mathematics.
\end{abstract}

\begin{center} {\bf FRACTALS, Vol 17 (1), 45-52, (2009)} \end{center}

\newpage

\section{Introduction}

Fractal sets are generally subsets of an Euclidean space with some fine
structures. These fine structures are usually indicative of nonsmooth
properties in the set. A prototype of such a set (in the one dimensional space
R of real numbers, for instance) is a Cantor set $C$ which is a compact and  totally disconnected set ( in the usual topology inherited from $R$) with vanishing topological dimension, although Card($C)=c$, the continuum. The Lebesgue measure of the set is also zero. Methods of ordinary analysis break down (or at the most yield unintuitive results) on such a set. Consider the function $s\left(  x\right)  =xf_{c}\left(
x\right)  $ where $f_{c}$: $\left[  0,1\right]  \rightarrow\left[  0,1\right]
$ is the Cantor's function. Then $\frac{d}{dx}\left(  \frac{s\left(  x\right)  }%
{x}\right)  =0$ a.e., although the function $\frac{s\left(  x\right)
}{x}$ is, nevertheless, non-constant as $\inf$ $f_{c}\left(  x\right) =0$, but $\sup f_{c}\left(x\right)=1$. The function $s\left(  x\right)  $ cannot therefore be
considerd to be a solution of an ordinary differential equation. Our aim in the work
is to give a reinterpretatin of the variability of the  Cantor function in the
framework of a scale invariant real analysis.

Recently there has been some interests in developing a framework of analysis on
a fractal set [1-3]. Parvate and Gangal [3], for instance, considered the so called  staircase functions, having a Cantor function like properties, in their formulation of the analysis. Their approach is based mainly on developing a formalism for replacing the linear Lebesgue measure (variable) viz.,  $x\in C \subset\left[  0,1\right]  $ by a nonlinear Hausdorff measure theoretic variable, viz., the integral staircase function $S_{c}^{s}\left(  x\right)  \approx x^{s}$ when $x (\approx 0)\in C$ and $s$ is the Hausdorff dimension of $C$. In the following, we outline the framework of a scale invariant analysis [4,5] based on a non-archimedean valuation [7]. We present mainly the measure theoretic aspects of the analysis. Besides identifying the Cantor function as a universal solution of a scale invariant first order ODE in the non-archimedean sense, we also show how the linear measure of $R$ gets extended naturally to the Hausdorff measure via a multiplicative model of $R$, equipped with the non-archimedean valuation.

In Section 2, we define a non-archimedean absolute value and corresponding valued measure on a zero (Lebesgue) measure Cantor set. The concept of {\em valued} (relative) infinitesimals [8] is introduced, which are then used to define the nontrivial valuation on $C$.  The valued measure is shown to be equal to the finite Hausdorff measure of the set. In Section 3, we present two examples on the triadic Cantor set and explain various properties of valued infinitesimals and related concepts. The nontrivial valuation is shown to be related to the Cantor function. In Section 4, we briefly outline the framework of  the scale invariant real analysis, and give simple examples of differentiation and the theory of Riemann integration. The validity of the mean value theorem is established in this non-archimedean setting and the Cantor function is realised as a {\em locally constant} function. The standard real integrals are also shown to have nontrivial asymptotic corrections from the scale invariant valued infinitesimals which might get significant applications in number theory and related subjects.

\section{Valuation and Measure}

A Cantor set $C$ is the unique limit set of an iterated system of similitudes
$f=\{  f_{i}|i=1,2,..,p\} \ :I\longrightarrow I,\  I=[0,1]$. For the triadic Cantor
set \ $p=2$ and we have
\begin{equation}
f_{1}\left(  x\right)  =\frac{1}{3}x,\ f_{2}\left(  x\right)  =\frac{1}%
{3}\left(  x+2\right)  ,x\in I\ \label{eqn1}%
\end{equation}
By definition, $C=f\left(  C\right).$ For simplicity, we assume in this paper that, the
similitude $f$ divides, for instance, $I$ into p  closed intervals each of
length ${1}/{r},$ and delete $q$ number of open intervals, so that
$p+q=r$. The scaling ratio of the similitude thus equals   ${1}/{r}$. More
complicated forms of $C$ will be considered separately. The Hausdorff
dimension of $C$ then turns out to be $s=\frac{\log p}{\log r}= \log_r p$.
The $s$ dimensional Hausdorff measure of $C$ \ equals 1, viz: $H^{s}\left(  C\right)  =1.$ The scaling law of Hausdorff measure also yeilds 
\begin{equation}\label{dim}
p=r^{s}
\end{equation}

Let $x\in C\subset I\subset R$. Then $x$ inherits the ordinary Euclidean
$\left(  \text{Archimedean}\right)  $ valuation $\left(  \text{absolute
value}\right)  $ $\ v:C\rightarrow R_{+}$ such that $v\left(  x\right)
=|x|=x$. Identifying $x\in C$ with the open set $\left(  0,x\right)\subset I, \ v\left(
x\right)$ may be considered to be the Lebesgue measure $v\left(  x\right)  =\mu
_{l}[\left( 0, x\right)]  $ with support $I\backslash C$. Since $I\backslash C+C=I$
and $\mu_{l}\left(  I\backslash C\right) = v(1) =1$, it follows that $\mu_{l}\left(
C\right)  =0$. The Hausdorff $s-$measure, on the other hand,  assigns $C$ a non trivial
measure, in the sense of a `content' or a uniform mass function 
\begin{equation}\label{smeas}
\mu_{s}[C] =  \lim_{\delta \rightarrow 0}\inf \sum_i (d(U_i))^s =1
\end{equation}
where $d(U)$ is the diameter of the set $U$ and the infimum is taken over all countable $\delta-$ covers $I_{\delta}=\{U_i\}$ such that $C\subset \bigcup U_i$, for the unique value of  $s$ given by equation (\ref{dim}). Although, both the Lebesgue measure and the s measure of a singleton $\{x\},\ x\in C$ is zero, the above scaling equation tells that on every application of the iterated system $f$, the total content (say, 1) of the set $C_n$ at the $n$th level of iteration is distributed uniformly over $p$ equal fragmented sets, each having a value $1/r^s$ so that $p\times(1/r)^s=1$. The total content of the original set $I$ therefore remains invariant for the final limit set $C=\bigcap C_n$. As a consequence, {\em the Hausdorff's $s-$measure gives, in a sense,  an idea of a relative measure: how the total content of the set gets fragmented over two consecutive iterated sets}, and thus provides, a more {\em intrinsic} sense of measure, in contrast to the Lebesgue  (outer) measure, which is an {\em extrinsic} concept, since the Euclidean norm is defined in relation to an exterior  reference point 0, say.  But for $s=1$, both these measures coincide, and the metric considered in $s$ measures is also equivalent to the Euclidean metric.  In the following we give a construction to reinterprete the Hausdorff measure purely in an intrinsic sense.

To this end, let us  now associate with x a new non-archimedean valuation as follows: 

{\bf Definition 1:} The disconnectedness of $C$ tells that to each $x\in C$, $\exists \ I_{\epsilon}(x)= \left(x-\epsilon,x+\epsilon\right)\subset I, \epsilon>0  $ (the largest open interval containing $x$) such that $\ C\cap I_{\epsilon}(x) =\{x\}$. Points in $I_{\varepsilon}\left(  x\right)  $ are said to be the {\em (relative) infinitesimal neighbours} in $I$ of $x\in C$. For an arbitarily small
$x$, on the other hand, $I_{\inf}^{+}=\underset{\varepsilon \neq 0}{\cup}I_{\varepsilon}\left(0\right)  =\underset{\varepsilon\neq 0}{\cup}\left(  0,\varepsilon\right)  $ is the open set of {\em (relative) infinitesimals}. For definiteness, $\varepsilon$ is said to denote a {\em scale}.

{\bf Definition 2:} Given an arbitrarily small  $x\in C$, there exists a suitable $\varepsilon$ and a set of relative infinitesimals $\tilde x\in I\backslash C$, such that $0 < \tilde x < \varepsilon< x,\ \tilde x= \lambda\varepsilon, \ ( 0 <\lambda< 1$) so that $\varepsilon/\tilde x\propto x/\varepsilon $ (in $I$).  Define the mapping 
(in I)
 ${v}:I_{\inf}^{+}\rightarrow (0,1)$ by

\begin{equation}
{v}\left(  \tilde x\right)  =\log_{\varepsilon^{-1}} 
{\frac{\varepsilon}{\tilde x}}, \label{value1}
\end{equation}
As a result, one obtains $\tilde x=\varepsilon^{1+v(\tilde x)},\ \varepsilon\rightarrow 0^+$. 

{\bf Remark 1}: $\varepsilon$, in general, may be the nearest rational approximation of a (positive) real $x \ (\ll 1)$. The {\em infinitesimals} are reals  $\tilde x < \varepsilon$ undetectable at the prescribed accuracy level determined by $\varepsilon$ and satisfy, by definition, the proportionality {\em (the law of inversion)} $\varepsilon/\tilde x\propto x/\varepsilon $. The definition (\ref{value1}) could further be justified as follows. We have, for any $ \tilde x=\varepsilon(1-\eta/\varepsilon),\ \eta=\varepsilon -  \tilde x$. For $\varepsilon\rightarrow  \tilde x$, which, in turn, $\rightarrow 0$, the ratio $\eta/\varepsilon$ becomes indeterminate and could assume a finite nonzero value. For a nonzero $\tilde x$, however, this vanishes.

{\bf Lemma 1:} The mapping ${v}$ is a non-archimedean valuation in the vector space of infinitesimals $I_{\inf}^{+}$ over $R$ ( considered as non-archimedean relative to the trivial valuation ${v}(r)=1, \ \forall r\in R$).

To prove this, note, first of all, that the relative infinitesimals $\overset{\sim}x$ corresponding to a given  $\varepsilon$  satisfy (by definition) the scale invariant equation
\begin{equation}\label{sfe}
x{\frac{d\overset{\sim}x}{dx}}=-\overset{\sim}x
\end{equation}
The corresponding solution space is a vector space on $R$.

Next, we show that, $0\neq\alpha\in I_{\rm inf}^+\backslash {C} \Rightarrow{v}\left(
\alpha\right)  =0$. Indeed, for such an $\alpha$, there exists $\varepsilon$ arbitrarily close to $\alpha$ so that ${v}\left(  \alpha\right)  =\underset{\varepsilon \rightarrow \alpha}{\lim}\log_{\varepsilon^{-1}}{\frac{\alpha}{\varepsilon}}=0$, by Remark 1. This is consistent with the fact that $\varepsilon=0$ is the only infinitesimal in $I$ relative to itself. It also follows that the definition may be extended to include $v(0)=0$, by setting $\varepsilon=0$ when $\alpha=0$.

Finally, we have for $x\in I_{\inf}^+ \cap C$ and a given scale $\varepsilon$

({ a}) ${v}\left(  x\right)  > 0, \  0\neq x\neq\varepsilon$. Further,  ${v}\left(  x\right)$ is also well defined.

(b) ${v}\left(  \alpha x\right)  ={v}\left(x\right)  ={v}\left(  \alpha\right) {v}\left(
x\right)$, since $v(\alpha)=1$ for an  $\alpha \in R$.

(c) To prove the strong triangle inequality $v(x_1+x_2)\leq \max \{v(x_1),v(x_2)\}$, let us  consider two infinitesimals $0<x_1\leq x_2<\varepsilon$, sufficiently small so that $0<x_1\leq x_2<x_1+x_2<\varepsilon$, for a preassigned scale $\varepsilon$. We note that this is a necessary condition for relative infinitesimals. We then have $v(x_1)\geq v(x_2)$ and so $v(x_1+x_2)= \log_{\varepsilon^{-1}} \varepsilon/(x_1+x_2) \leq v(x_1)$.

The set of relative infinitesimals $I_{\inf}$ is thus an ultrametric vector space [7],
with $v$ being a nonarchimedean absolute value over infinitesimals.

{\bf Remark 2:} Even as $x\rightarrow 0$ in $C$, the valuation $v(x) = \log_{\epsilon^{-1}}\lambda^{-1}$ assigns a unique value to a one parameter family of infinitesimals relative to $\varepsilon$ ( notice the proportionality constant in the 
inversion law) to  $x$. However, the scale $\varepsilon$ may be any number in (0,1). Consequently, the set of (relative) infinitesimals is of cardinality $c$, the continuum. 

To examine the structure of $I_{\inf }$, let us recall the basic topological properties of a non-archimedean space.

(i) Every open ball $B_r(a)=\{x| v(x-a)<r\}$ is closed and vice versa. Such a ball is called a clopen set.

(ii) Every point $b\in B_r(a)$ is a centre of $B_r(a)$.

(iii)Any two balls in $I_{\inf }$ are either disjoint or one is contained in  another.

(iv) Any open set, and hence $I_{\inf }$ itself, is the union  at most of a countable family of closed balls.

The above results are simple consequences of the ultrametric inequality (c). The property (iv), in particular, holds because $I_{\inf }\subset R$ is the union of at most countable number of open balls, which are nevertheless closed in the present context.  It also follows that $I_{\inf }$ with the above norm is a totally disconnected set. Moreover, every closed ball is compact (c.f. Lemma 2) and so is covered by a finite number of closed balls in each of which $v(x)$ is a constant (Remark 2). Given a scale $\varepsilon$, these finite partition introduces a further set of finer scales, $\varepsilon^j\ j (\geq 1)\in J$, a finite subset of $R$. 

{\bf Lemma 2:} A closed ball in $I_{\inf}$ is both complete and compact.

The proof follows from the following observations. Given an $\varepsilon >0$, consider a closed interval $[a,b] \subset I_{\inf}$ (in the usual topology)  such that $0<a<b<\varepsilon$. The valuation $v$ maps this closed interval onto the closed interval (ball) $B_r(c)=[\tilde a,\tilde b]\subset I_{\inf}$ where $\tilde a=v(b)$ and $\tilde b =v(a)$, with centre $c=(\tilde a+\tilde b)/2$ and radius $r=v(\tilde b-\tilde a)$. Completeness follows from the standard ultrametric properties: $\{x_n\}$ is Cauchy $\Leftrightarrow v(x_n-x_m)\rightarrow 0 \ \Leftrightarrow v(x_n-x_{n-1})\rightarrow 0 \ \Leftrightarrow x_n=x_N$ for $n\geq N$. Compactness is also verified directly following the standard metric space techniques.

{\bf Remark 3:} The ultrametricity tells that $v(x)=\sigma^{ a(x)}, 0< \sigma <1$ for a suitable $\sigma=\sigma(\varepsilon)$ and $a(x)$ satisfying $a(x+y)\geq \min \{a(x),a(y)\}$.  For definiteness, (and without any loss of generality) let us fix $\sigma$ by choosing $\sigma(\varepsilon)=\varepsilon$,  so  that $v(x)=\varepsilon^{a(x)}, 0<\varepsilon <1$. In relation to  the finite disjoint cover, and given the choice of the scale $\varepsilon: 0<\varepsilon<1$, an infinitesimal $x$ is now parametrised as $x=x_n\tilde x$ where $x_n=\varepsilon ^{1+\varepsilon^{ a(x_n)}}, \ 0<x<\varepsilon $ and $\tilde x=1 + X(x)$ so that $v(x)=v(x_n)= \log_{\varepsilon ^{-1}}x_n/\varepsilon$ and $v(\tilde x)=1$, since $\log_{\varepsilon^{-1}}\tilde x=0$ as $\varepsilon\rightarrow 0$, for a suitable $X(x)$. Also $a(x_n)= a_n,$ for a sequence of constants $a_n $ (by Remark 2), from a finite set. Accordingly, infinitesimals in the $n$th covering ball have the constant valuation $v(x_n)=v_n$, say. For latter convenience, we make these $n$ and $\varepsilon$ dependences of $v$ explicit by writing $v_n=\alpha_n\varepsilon^{s_0}$ so that $\alpha_n$ now assumes values from the finite set and $s_0$ is another constant.  Moreover, given an $\varepsilon$, one can always choose an infinitesimal $x$ so that $v(x)\leq x$. Stated in other words, given an infinitesimal $x>0$, we can choose a scale $\varepsilon$ so that $v(x)\leq x$. With this assumption, for $[a,b]$, we have $v(b-a)\leq |b-a|$.  More details of the valuation will be given in Example 1 in the next section. Summing up the above observations, we now have

{\bf Proposition 1:} The mapping $v: I_{\inf}\rightarrow (0,1)$ is a nonarchimedean valuation  on the set $I_{\inf}$ and assumes values from a finite set of (0,1).  It has the general form $v(x)=\alpha_n\varepsilon^{s_0}, \ x\in I_{\inf}$ for a given choice of the scale $\varepsilon$, $\alpha_n$ and $s_0$ being two constants.

The proof follows from the representation of infinitesimals in the above Remark 3.

We now utilize this valued set of infinitesimals to define a valued measure on
$C$ in several steps.

$\left(  i\right)$ Each $(0\neq) \ x\in C$ is identified with $I_{\inf}\left(
x\right)  =x+I_{\inf},\ I_{\inf}=I_{\inf}^{+}\cup I_{\inf}^{-}, \ I_{\inf}%
^{-}=\underset{\epsilon}{\cup}\left(  -\epsilon,0\right)  $.
Clearly, $\underset{x}{\cup}I_{\inf}\left(  x\right)  \supset I$.

$\left(  ii\right)$  Given $x\in C$, define a (one parameter family of ) multiplicative neighbour(s) in $I_{\inf}\left(  x\right)  $ which are induced by the valued infinitesimals $\tau$ in $I_{\inf }$ by

\begin{equation}
I_{\inf}\left(  x\right)  \ni X_\tau^{\pm}=x.x^{\pm v(\tau)}=x.x^{\pm \alpha_n \varepsilon^s}\label{multnebor}
\end{equation}
where $\alpha_n=\alpha_n(x)$ may now depend on $x$.
This is a generalization of the trivial equality $x=x.x^{0}$ in $I$.

$\left(  iii\right)$ Define a new absolute  value of $x\in C$ by

\begin{equation}
||x||=\inf \log_{x^{-1}} {X_+}/{x}=\inf \log_{x^{-1}}{x}/{X^{-}}\label{value2}
\end{equation}
so that $||x||=\varepsilon^{s}$, where $\varepsilon^s=\inf \alpha_n\varepsilon^{s_0}$ and the infimum is over all n. The absolute value thus picks the maximum of the finer scale neighbours $X_\tau^ {\pm}$. As shown below, the number $s$ equals the Hausdorff dimension of the set. It now follows also that 

{\bf Proposition 2:} $||.||:C\rightarrow R_{+}$ is a non-archimedean valuation.

{\bf Remark 4}: Both $I_{\inf}$ and $I_{\inf}\left(  x\right)  $ are totally
disconnected in the topology induced by the non-archimedean valuation (\ref{value2}), although both are connected as subsets of $R$ in the usual topology. Given a totally disconnected set $C$ (in usual topology), the non-archimedean topology induced by $v$  injects a finer subdivisions into the infinitesimal neighbourhood of $0^+$, which is then inherited by neighbourhoods of any finite $x\in C$. However, the $||.||$-topology misses these local fine structures, since the $||.||$ metric reads off only the supremum of $I_{\inf}(x)$. However, as $x \rightarrow 0^+$, the two topologies may coincide, provided $n=1$ and $\alpha_1=1$.

We now define the valued measure $\mu_{v}: \ C\rightarrow R_{+}$ by

(a) $\mu_v\left(  \Phi\right)  =0, \ \Phi$ the null set.

(b) $\mu_{v}[\left(  0,x\right)  ]=||x||$, when $x\in C$.

(c) For any $E\subset C$, on the other hand, we have $\mu_v(E)= \underset{\delta \rightarrow0}{\lim}\inf \ \sum_i  \{d_{\rm na}(I_i)\}$, where $I_i\in \tilde I_{\delta}$ and the infimum is over all countable $\delta$- covers $\tilde I_{\delta}$ of $E$ by clopen balls. Moreover, $ d_{\rm na}(I_i)$=the non-archimedean  diameter of $I_i= \sup\{||x-y||:\  x, \ y \ \in I_i\}$. Denoting the diameter in the usual (Euclidean) sense by $d(I_i)$, one notes that $d_{\rm na}(I_i)\leq \{d(I_i)\}^{s}$, since  $x,\ y\in C$ and $|x-y|= d$, imply  $||x-y||= \varepsilon^{s}\leq d^s$, as the scale $\varepsilon$ satisfies, by definition, $\varepsilon\leq d\leq\delta$.

Thus $\mu_v$ is a metric (Lebesgue outer) measure on $C$ realised as a non-archimedean ultrametric space. Now to compare this with the Hausdorff $s$ measure, we first note that $\mu_v[E] \leq  \mu_{s}[E]$ since  $d_{\rm na}(I_i)\leq \{d(I_i)\}^{s}$ for a given cover of (Euclidean) size $\epsilon$. Next, for a cover of  clopen balls of sizes $\epsilon_i$, we have  $\sum_i  \{d_{\rm na}(I_i)\} =\sum_i  \epsilon^s_i$. For the Hausdorff measure, on the other hand, covers by any arbitrary sets are considered. Using the monotonicity of measures it follows that  
 
\begin{equation}\label{meas1}
   \inf \ \sum_i  \{(d(I_i))^{s}\}\leq \inf \ \sum_i  \{d_{\rm na}(I_i)\}
 \end{equation}
so that letting $\epsilon\rightarrow 0$ we have $\mu_v[E]\geq \mu_{s}[E]$. Hence 

\begin{equation}\label{meas2}
 \mu_v[E]= \mu_{s}[E]
 \end{equation}
for any subset $E$ of $C$. Finally, for $s=$ dimension of $C$, $\mu_s[C]$ is finite and hence the valued measure of $C$ is also finite. Notice that the valued measure  selects {\em naturally} the dimension of the Cantor set. 

\section{Examples}

{\bf 1:} Let us now investigate in detail the well known triadic Cantor set $C$ in the light of the above analysis. Suppose we begin with the  set $C_0=[0,1]$. In relation to the {\em scale} 1, $C_0$ is essentially considered to be a doublet $\{0,1\}$, in the sense that real numbers $0<x<1$ are {\em undetectable} in the assigned scale, and  hence all such numbers might be identified with 0. We denote this 0 as $0_0= [0,1)$, the set of {\em infinitesimals}. However, the possible existence of {\em infinitesimals} are ignored at this scale and so 0 is considered simply as a singleton $\{0\}$ only.
 At the next level, we choose a smaller scale $\varepsilon=1/3$ (say), so that only the elements in $[0,1/3)\subset C_0$ are now  identified with 0, so that $0_1=[0,1/3)$, which is actually $0_1=0_0$ in the unit of 1/3. Relative to this nontrivial scale 1/3, we now assign the ultrametric valuation $v$ to $0_1$. In principle, all possible ultrametric valuations are admissible here. One has to make {\em a priori choice} to {\em select } the most appropriate valuation in a given application. In the context of the triadic Cantor set, there happens to be a unique choice relating it to the Cantor's function, as explained below.

Recall that the valuation induces a nontrivial topology in $0_1$. Accordingly, the set is covered by $n$ number of disjoint clopen intervals of {\em valued} infinitesimals. At the level 1, $n=1$, which is actually the clopen interval $I_{11}$ of length 1/3 and displaced appropriately to the middle of the 1/3rd Cantor set, viz, $I_{11}=[1/3,2/3]$ (in the ordinary representation this is the deleted open interval, including the two end points of neighbouring closed intervals). The value assigned to these valued infinitesimals is the constant $v(I_{11})=1/2$, where, of course, $v(0)=0$. In principle, again, $v$ could assume any constant value. Our choice is guided by the triadic  Cantor function. Thus the valued set of infinitesimals, at the scale 1/3, turns out to be $0_1=\{0,1/2\}$. 

How does this valued set of infinitesimals enlight the ordinary construction of the Cantor set? Let $C_1= F_{11}\cup F_{12}$ where $F_{11}=[0,1/3]$ and $ F_{12}=[2/3,1]$. The value awarded to the deleted middle open interval is now inherited by these two closed (clopen) intervals, and so $||F_{11}||=1/3^s$ and $|| F_{12}||=1/3^s$, recalling that $2=3^s$, $s$ being the Hausdorff dimension $s=\log2/\log 3$. 

At the next level, when the scale is $\varepsilon=1/3^2$, the above interpretation can be easily extended. The zero set is now made of 3 clopen sets $0_2= I_{20}\cup I_{21}\cup I_{22}$ where $I_{20}=[1/9,2/9],\ I_{21}= [3/9,6/9]$ and $I_{22}=[7/9,8/9]$. The value assigned to each of these sets are respectively, $v(I_{20})=1/4,\ v(I_{21})=2/4$ and $v(I_{22})=3/4$, so that the valued infinitesimals are given by $0_2=\{0, 1/4,2/4,3/4\}$. Notice that the new members of the valued infinitesimals are  derived as the mean value of two consecutive values from those (including 1 as well) at the previous level. These valued infinitesimals now, in turn, assign equal value to the 4 closed intervals in the ordinary level 2 Cantor set $C_2 = F_{20}\cup F_{21}\cup F_{22}\cup F_{23}$ where $F_{20}=[0,1/9]$ and etc, viz. $||F_{2i}||=1/2^2=1/3^{2s}, \ i=0,1,2,3$. Notice that, in the sense of Sec. 2, the valued infinitesimals $0_2$ induces a {\em fine structure} in the neighbourhood of $F_{2i}:$ for a $x\in F_{2i}$, we now have valued neighbours $X^{\pm}=xx^{\pm k3^{-2s}}, \ k=1,2,3$. Clearly, $||F_{2i}||=||x||=1/3^{2s}$, the infimum of all possible valued members, so misses the above fine structures (c.f., Remark 4). It also follows that the limit set of this triadic construction  reproduces the Cantor function (c.f., Example 2) as the the valuation $v: [0,1]\rightarrow [0,1]$,  defined originally on the {\em inverted} Cantor set ${\bf 0}=\bigcap_{n}\bigcup_{k}I_{nk}$, and then extended on [0,1] by continuity.

{\bf Remark 5:} The continuity in the present ultrametric topology is defined in the usual manner. Further, $v$ on $\bf 0$ is an example of {\em locally constant} function relative to the $||.||$-topology and will be shown (in the next section) to satisfy the differential equation  

\begin{equation}\label{de1}
x{\frac{d v(x)}{dx}}=0
\end{equation}
We may interprete this as follows: Considered as a function on ${\bf 0}$ (or $C$), $v$ is constant in clopen sets $I_{nk}( \ \rm{or} \ F_{nk})$ for fixed values of both $n$ and $k$, but experiences variability as either of these vary. This variability is not only continuous, but continuously first order differentiable as well. In constrast, $v$ on $\{I_{nk} \ {\rm or} \ ( F_{nk})\}$ is a discontinuous function in the usual topology.

 {\bf 2:} In this example, we present an explicit construction of multiplicative neighbours of $x\in C$ using the Cantor function $f_C:I\rightarrow I$. In the following we denote this function instead by $\tilde X(x)$.
Consider the $\frac{1}{3}$-rd Cantor set: $r=3,\ p=2.$
Let $x={\Sigma}{a_{i}}{3^{-i}}$ be the ternary repesentation of
$x\in C$ where $a_{i}$  may be either 0 or 2. We set $x=\frac{2}
{3}\psi(  \tilde{X})  $ where $\psi(  \tilde{X})  ={\Sigma}\frac{b_{i}}{3^{i-1}}$ and $\tilde {X}=\Sigma b_{i}2^{-i}\in I\backslash C, \ b_{i}\in\{0,1\}$.

Then $\tilde{X}=\tilde{X}\left(  x\right)  $ defined as the
inverse of the above functional equation is the Cantor function
$\tilde X : [0,1]\rightarrow\lbrack0,1]$. By continuity, this extends over $C$ as well.

Let us recall that at the k -th step of the iterative construction
of the Cantor set, the initial closed interval $I$ fragments into $2^{k}$
smaller closed intervals $I^{k}_{j}= [x_{2j-1},x_{2j}], \ j=1,2....,$ each of length $3^{-k}$. 

Then $x_{2j}-x_{2j-1}=3^{-k}$. Definition of the Cantor function also gives that 

\begin{equation}
\tilde{X}{(x_{2j})-\tilde X(x_{2j-1})=2^{-k}}\label{eqna}
\end{equation}

Let $x\in C$. Then $x\in I_{j}^{k}$ for some j. It thus follows
\begin{equation}
\tilde{X}(x_{2j})-\tilde X(x_{2j-1})=\frac{3^{k}}{2^{k}}(x_{2j}%
-x_{2j-1})\label{eqnb}%
\end{equation}

Let $\tilde X(x_{2j})=X_{+}, \ \tilde X(x_{2j-1})=X_{-}, \ x_{2j}=x_{+}, \ x_{2j-1}=x_{-}$.
Suppose also that

$$
{3^{k}(x_{+}-x)\rightarrow k\log\sigma_{+},\ \ 3^{k}(x-x_{-})\rightarrow
k\log\sigma_{-}} 
$$
and
\begin{equation} 
2^{k}(\tilde{X}_{+}-\tilde{X})\rightarrow k\log X_{+}
^{^{\prime}},\ 2^{k}(\tilde{X}-\tilde{X_{-}})\rightarrow k\log X_-^{\prime}
\end{equation}
for infinitely large $k \rightarrow\infty$. The limiting value of (\ref{eqnb})
thus becomes%

\begin{equation}
\log X_{+}^{\prime}+\log X_{-}^{\prime}=\log\sigma_{+}+\log\sigma
_{-}\label{eqnc}%
\end{equation}

Now, using the inequality $\frac{\alpha+\gamma}{\beta+\delta}\leq\max(\frac{\alpha
}{\beta},\frac{\gamma}{\delta}),\alpha,\gamma\geq0,\beta,\delta>0$, 
(\ref{eqnc}) yields
\begin{equation}
\max\left(  \frac{\log X_{+}^{^{\prime}}}{\log\sigma_{+}},\frac{\log
X_{-}^{^{\prime}}}{\log\sigma_{-}}\right)  \geq1\label{eqnd}%
\end{equation}

But $\left(  14\right)  $ shows that $\sigma_{+}=\sigma_{-}^{-1}=\sigma$ (say) and 
$X_{+}=X_{-}^{-1}$, as $k \rightarrow \infty$, so that (15) reduces to

\begin{equation}
X_{+}^{^{\prime}}=\sigma^{1+j},\ X_{-}^{^{\prime}}=\sigma^{-(1+j)}
, \ j\geq0\label{eqne}
\end{equation}

Setting $\sigma^{-1}X_{+}^{^{\prime}}=\frac{X_{+}}{x}$ and $\sigma
X_{-}^{^{\prime}}=\frac{X_{-}}{x}$, 
we finally get the multiplicative neighbours of $x\in C$ as

\begin{equation}
X_{\pm}=x\sigma^{\pm j}\label{eqnf}
\end{equation}
Notice that $\sigma \approx 1$. In the notation of Section 2, $\sigma=x^{\tau^s}, \ \tau$ being a valued infinitesimal. Although the inequality eq(\ref{eqnd}) is reminiscient of the non-archimedean valuation, we are at this stage unable to restrict the set of values of $j$  directly from the above limiting argument.

\section{Differentiability}

The framework of elementary Calculus can now be developed on a Cantor set $C$ when the ordinary Cantor set $C$ is replaced by the valued Cantor set $\bf C$. Each element  of such a set is assigned a non-archimedean valuation $||.||$. A valued point $X\in \bf C$ now can undergo changes in the set {\em continuously} form one site $X_1$ to another $X_2$ by infinitesimal {\em hoppings} following {\em the law of inversion}, viz, $||X_2||=||X_1||^{\alpha}, \ 0<\alpha<1 $ ($X_2\geq X_1$, in the Euclidean sense). Continuity is defined in the standard manner using the metric induced by $||.||$. As an example, let us define the differentiability as follows: A mapping $f: \ C\rightarrow C $ is said to be differentiable at $X_0$ if $\exists$ a finite $l$ such that $0<||X-X_0||<\delta \ \Rightarrow \ | \ ||f(X)-f(X_0)||/||X-X_0||-l| < \epsilon$ for $\epsilon >0$ and $\delta(\epsilon)>0$ and we write $f^{\prime}(X_0)=l$.  Since the evaluation of $||.||$ valuation amounts to evaluting the infinitesimal valuation $v$ for an infinitesimal $x$ living in closed intervals of $I\subset R$, the mean value theorem of the form 
\begin{equation}\label{mvt}
f(X)= f(X_0)+f^{\prime}(X_0)||X-X_0|| +O(||X-X_0||^2) 
\end{equation}
is valid. Accordingly, it follows from eqns(\ref{value2}) and (\ref{value1}) that 
\begin{equation}\label{diff}
f^{\prime}(X)={\frac{ d\log  f(x)}{d\log x}}
\end{equation}
for a positive (valued) infinitesimal variable $x$. The scale invariance of the above is a consequence of the logarithmic derivative in $R$.  

We now show that the Cantor function, which in our notation is the valuation $v: I\rightarrow I, \ I= [0,1]$, is a locally constant function. Let us recall that a totally disconnected Cantor set is represented recursively as $C=\cap_n\cup_k F_{nk}$ (c.f. Sec.3). The set $I$ thus is written as $I =\cap_n\cup_k(F_{nk}^{\prime}\cup I_{nk})$, the open interval $F_{nk}^{\prime}$ being $F_{nk}$ with end points removed (recall that $I_{nk}$ are closed intervals). By definition, $v(I_{nk})=a_{nk}$, a constant for each $n$ and $k$. We set $v(F_{nk}^{\prime})=0$. This equality is to be understood in the following sense. For instance, in the triadic Cantor set, the zero value of, say, $v(F_{11}^{\prime})=0$ becomes valued finitely at the next level, viz, $v(I_{21})=1/4$ but $v(F_{21}^{\prime})=v(F_{22}^{\prime})=0$. The derivative of $v$ vanishes not only for each $n$ and $k$, but also as $n\rightarrow \infty$, because every Cantor point (point of discontinuity, in the ordinary sense) is replaced by the open interval (-1,1) of valued infinitesimal neighbours, at an appropriate scale $\varepsilon$, which persists even as $\varepsilon\rightarrow 0$. Stated in other words, the above forms of block decompositions of $I$ into infinitesimally small, valued scales, are available even as $n\rightarrow \infty$. This explains the removal of the zero measure discontinuous set of the derivative of the Cantor function in the present scale invariant analysis. We remark that to catch the finer scale variability of $v$, one needs to work with the topology induced by $v$ itself on the set of infinitesimals $I_{\inf}$. Extending the above definition of differentiability in the $v$ metric, one may verify that $dv/d\log x^{-1}=1$ for $x\in I_{\inf}$, instead.

Let us mention also, for the sake of completion, that integrals of continuous functions could be defined analogously. We remark that the above structure of a scale invariant calculus could also be developed even in $R$ in the presence of nontrivial valued infinitesimals. Indeed, the relevant definitions could be extended  in the following manner. One may replace the singleton set $\{0\}$ of $R$ by {\em a} Cantor set $C$. A set of scales is thus generated, relative to which one can then introduce a set of valued infinitesimals in $R$ for each choice of the Cantor set $C$. The real number set $R$, in this new representation, accquires the structure of a measure one Cantor set. The determination of the precise structure of the Cantor set will be considered separately. In Ref.[4] (and in references cited there), a fractal dimension of the extended real number system $\bf R$ is estimated to be 1+$\nu$, where $\nu=\frac{\sqrt 5-1}{2}$, the golden ratio, under the assumption that the increments in $\bf R$ is accomplished as an SL(2, R) group action. As an example, we evaluate  $\int^1_0 dx$ in such an extended $\bf R$.  We have
\begin{equation}\label{int}
 \int^1_0 dx= \underset{\epsilon \rightarrow 0}{\lim}\left[\int^1_{\epsilon} + \int^{\epsilon}_{\epsilon^{1+v(\epsilon)}} dx \right]= \underset{\epsilon \rightarrow 0}{\lim}\left[1-\epsilon +v(\epsilon)\right]\approx 1+v(\epsilon)
\end{equation}
when scales less than $\epsilon$ are replaced by the valued scales $\epsilon^ {1+v(\epsilon)}$ and arbitrarily small $x$ is replaced by the corresponding nontrivial value $v(x)=\log_{\epsilon^{-1}}(\epsilon/x)$. The above form of  the well known integral corresponds to an asymptotic correction, which might become significant in  analysis and number theory, for instance, in the prime number theorem [6]. 

We note finally that the nondifferentiability of $|x|$ at 0, in the usual topology,   gets {\em smoothed out} in the $||.||$ topology, since the ordinary zero is replaced by valued infinitesimals of the form $\varepsilon^s$ relative to the scale $\varepsilon$, so that  both right and left hand (logarithmic) derivatives of eq(\ref{diff}) equal 1, as in the case for any finite $x$.


\newpage

\begin{center}
{\bf Earratum to ``Analysis on a fractal set''\\
{Fractals, vol 17, No. 4 (2009) 547} }
\end{center}

1. Page 46, 4th line, column 1: the function $s(x)=x_c \ \rightarrow xf_c(x)$ should be replaced by 
$s(x)=xf_c(x)$.

2. Page 46, Sec 2, 2nd line, column 2: $f=\{ f_i|i=12\ldots p\}$ should be replaced by 
$f=\{ f_i|i=1, 2, \ldots , p\}$.

3. Page 50, 5th line in Example 2( column 1), $r={\frac{1}{3}},\ p=2$ should be replaced by $r=3, \ p=2$.

{\bf An Explanatory Remark:} (Page 47, column 2)  In paragraph 2, following equation (5), we show that $0\neq \alpha \ \in I^+_{\inf}/C \ \Rightarrow v(\alpha)=0$. This might appear to conflict Definition 2 for the nontrivial valuation for a relative infinitesimal $\tilde x \ \in I/C $ for an arbitrarily small $x \ \in C$. However, there is a tacit assumption, suggesting that nontrivial infinitesimals $\tilde x$ live in a finite number of disjoint closed intervals of $(0, \varepsilon)$ and satisfy the inversion rule (4). Infinitesimals $\alpha$, on the otherhand, are {\em real like} in the sense that these do not respect the inversion rule. Further, the scale $\varepsilon$ could be translated arbitrarily close to $\alpha$. This assumption is justified, a posteriori, in the examples of Sec. 3, and also in our recent works [1,2].

{\bf References}

1. S Raut, D P Datta, Nonarchimedean scale invariance and Cantor sets, FRACTALS, March (2010), to appear.

2. D P Datta, A Roy Choudhuri, Scale free analysis and prime number theorem, Fractals, to appear.

\end{document}